\documentclass[10pt]{article}
\usepackage{times}

\usepackage{latexsym}
\usepackage{amsmath}
\usepackage{amssymb}

\newtheorem{Theor}{Theorem}[section]
\newtheorem{Prop}[Theor]{Proposition}

\newtheorem{Rem}[Theor]{Remark}
\newtheorem{Lemma}[Theor]{Lemma}
\newtheorem{Not}[Theor]{Notation}
\newtheorem{Def}[Theor]{Definition}

\newtheorem{Cor}[Theor]{Corollary}

\def\Dem{\hbox{\textsc {Proof.}\,\,}}
\def\Aut{\hbox{\sffamily {Aut}\,}} 
\def\Id{\hbox{\sffamily {Id}\,}}
\def\Tri{\hbox{\sffamily {Tri}\,}}
\def\Pro{\hbox{\sffamily {Pro}\,}}
\def\Min{\hbox{\sffamily {Min}\,}} 
\def\rank{\hbox{\sffamily {rank}\,}} 
\def\one{\hbox{\sffamily {1}\,}} 
\def\Sym{\hbox{\sffamily {Sym}\,}}
\def\Lin{\hbox{\sffamily {Span}\,}}

\def\qed{\hfill $\Box $}

\usepackage{version}

\newenvironment{mm}
  {\marginpar{$\Downarrow$}
}
  {\marginpar{$\Uparrow$}}

\def\nomarks{\def\mm{}\def\endmm{}}

\nomarks

\begin{document}

\title{Manifolds of algebraic elements in the algebra
  $\mathcal{L}(H)$ \\ of bounded linear operators.}

\author{Jos\'e M. Isidro \thanks{Supported by Ministerio de 
Educaci\'on y Cultura of Spain, Research Project PB 98-1371.}\\
Facultad de Matem\'aticas, \\
Universidad de Santiago,\\
Santiago de Compostela, Spain.\\
{\normalsize\tt jmisidro@zmat.usc.es}
\and
Michael Mackey\thanks{Supported by Ministerio de 
Educaci\'on y Cultura of Spain, Research Project PB 98-1371.}\\
Department of Mathematics,\\
University College Dublin,\\
Dublin 4, Ireland.\\
{\normalsize\tt michael.mackey@ucd.ie}
}

\date{October 12,  2001}

\maketitle

\begin{abstract}
  Given a complex Hilbert space $H$, we study the differential
  geometry of the manifold $\mathcal{M}$ of all projections in
  $V=\mathcal{L}(H)$. Using the algebraic structure of $V$, a
  torsionfree affine connection $\nabla$ (that is invariant under the
  group of automorphisms of $V$) is defined on every connected
  component $\mathfrak{M}$ of $\mathcal{M}$, which in this way becomes
  a \begin{mm}symmetric holomorphic \end{mm} manifold that consists of
  projections of the same rank $r$, ($0\leq r\leq \infty$). We prove
  that $\mathfrak{M}$ admits a Riemann structure if and only if
  $\mathfrak{M}$ consists of projections that have the same finite
  rank $r$ or the same finite corank, and in that case $\nabla$ is the
  Levi-Civita \begin{mm}and the K\"ahler\end{mm} connection of
  $\mathfrak{M}$.  Moreover, $\mathfrak{M}$ turns out to be a totally
  geodesic Riemann manifold whose geodesics and Riemann distance are
  computed.  \bigskip


\end{abstract}


\section{Preliminaries on
JB-algebras}
\subsection{Introduction.}

In \cite{HIR} Hirzebruch proved that the manifold of minimal
projections in a finite-dimensional simple formally real Jordan
algebra is a compact Riemann symmetric space of rank 1, and that any
such space arises in this way.  Later on, in \cite{NOM1} Nomura
established similar results for the manifold of minimal projections in
a topologically simple JH-algebra (a real Jordan-Hilbert algebra). The
results in \cite{CHIS}, \cite{ISI} and \cite{ISTA} lead to the idea
that the structure of a JBW-algebra $V$ might encode information about
the differential geometry of some manifolds naturally associated to it
\cite{KAUJ}. In particular, that the knowledge of the JBW-structure of
$V$ is sufficient to study the manifold of projections in $V$. Every
JBW-algebra can be decomposed into a sum of closed ideals $V=V_I\oplus
V_{II}\oplus V_{III}$ of types I, II, and III respectively, and for
our purpose it is not a hard restriction to assume that $V$ is
irreducible. JBW-algebras of type III are not well
understood. A typical example of a type II JBW-algebra is
$L^{\infty}[0,1]$ whose lattice of projections is modular. Namely,
it consists of characteristic functions of Lebesgue measurable
subsets of $[0,1]$ and form a discrete topological space since we
have $\Vert \chi_1 -\chi_2\Vert =1$ whenever $\chi_1\neq \chi_2$. 
Thus we have to consider JBW-algebras of type I and we shall
assume them to be factors, hence factors of type
$I_n$ for some cardinal number
$1\leq n\leq \infty$. Those of type $I_2$, called spin factors, are
Hilbert spaces (\cite{HAOL} th. 6.1.8) hence they are included in
the work of Nomura. Factors of type
$I_n$ with $3\leq n<\infty$ are certain spaces of matrices
(\cite{HAOL} th. 5.3.8) and so they are included in the work of
Hirzebruch. Thus essentially we have to consider the JBW-algebra
$V\colon =\mathcal{L}(H)_{\rm sa}$ of the selfadjoint operators on a
Hilbert space $H$ over some of fields $\mathbb{R},
\,\mathbb{C},\,\mathbb{H}$ (\cite{HAOL} th.  7.5.11). Here we make
such a study in a systematic manner without the use of any global
scalar product in $V$, in the complex case. With minor changes it
applies to the other  two fields. 

The set $\mathcal{M}$ of all projections in $V$ can be
identified with the set of all closed subspaces of $H$, which is a
Grassmann manifold $\mathbb{G}(H)$ in a classical way \cite{KAUS}. It
is known that $\mathbb{G}(H)$ has several connected components
$\mathfrak{M}$, each of which consists of projections $p$ in $V$ that
have a fixed rank $r$, $0\leq r\leq \infty$. An affine connection
$\nabla$, that is invariant under the group $\Aut ^{\circ}(V)$ of
automorphisms of $V$, is then defined on each connected component
$\mathfrak{M}$ with only the help of the
JBW-structure. With it, $\mathfrak{M}$ becomes a symmetric totally
geodesic real analytic manifold. Moreover, it is possible (and in fact
easy) to integrate the differential equation of the geodesics
corresponding to initial conditions defined by purely algebraic
equations. For $r=1$, $\mathfrak{M}$ is the complex projective space 
$\mathbb{P}(H)$.

Motivated by the above, we ask whether it is possible to define a
Riemann structure on $\mathfrak{M}$, and a necessary and sufficient
condition for this to happen is established. The tangent space to
$\mathfrak{M}$ at a point $p$ is the range of the $\frac{1}{2}$-Peirce
projector of $V$ at $p$, and $\mathfrak{M}$ admits a Riemann structure
(if and) only if $P_{1/2}(p)V$ is (homeomorphic to) a Hilbert space.
In \cite{ISTA} it has been proved that the latter occurs if and only
if the rank of $p$ or the corank of $p$ (the rank of $\one -p$, where
$\one$ is the unit of the algebra $V$) is finite.  $\mathcal{M}$ is an
ortho-complemented lattice, and the mapping $p\mapsto p^{\perp}\colon
=\one -p$ is an involutory homeomorphism. In fact, this involution is
a real analytic diffeomorphism of $\mathcal{M}$, hence it suffices to
study the connected manifolds $\mathfrak{M}$ with $r< \infty$ which
leads us again to the work of Nomura. Namely, consider the algebras of
finite rank operators, of Hilbert-Schmidt operators, of compact
operators, and of bounded operators on $H$, respectively, and the
inclusions
$$
\mathcal{F}(H)\subset \mathcal{L}_2(H)\subset
\mathcal{L}(H)_0 
\subset \mathcal{L}(H)
$$
The first three of these algebras have the same set of projections,
which is exactly the set of {\sl finite rank} projections in
$\mathcal{L}(H)$. However, the topologies induced on $\mathcal{F}(H)$
by $\mathcal{L}_2(H)$ and $\mathcal{L}(H)$ do not coincide (unless
$\dim H<\infty$), and a priori there is no reason to expect that they
should coincide on the set of projections (we shall see that this
happens).

We then study the Riemann manifolds $\mathfrak{M}$ for $r<\infty$
without the help of any global scalar product in these algebras. A
scalar product in the tangent bundle to $\mathfrak{M}$ is needed, of
course, but it is locally provided in a canonical way by the
JBW-algebra structure of $V$. We begin 
with a discussion of the subalgebras $V[a, b]$ generated in $V$ by
certain pairs of elements $(a,\,b)$. These subalgebras, that play a
fundamental role in our study, turn out to be Jordan isomorphic to
$\Sym (\mathbb{R},\, 2)$, the algebra of $2\times 2$ symmetric
matrices with real entries and the usual Jordan matrix product, and
therefore they are finite dimensional. By choosing an appropriate
basis in $V[a,b]$ it is easy to integrate the differential equations
of the geodesics in $\mathfrak{M}$. 


\subsection{Preliminaries on JBW-algebras.} 

A Jordan algebra $V$ is an algebra over $\mathbb{R}$ or $\mathbb{C}$
in which the following two identities hold for all $x,\,y$ in $V$:
\begin{equation}\label{id1} 
xy=yx, \qquad  x^2(xy)=x(x^2y)
\end{equation} 
Let $V$ be a Jordan algebra.  Then $L(x)$ and $P(x)$, $(x\in V)$, are
defined by
\begin{equation}\label{id2}
L(x)y\colon =xy, \qquad 
P(x)y\colon = 2 L(x)^2y-L(x^2)y, \qquad (y\in V).
\end{equation}

An element $a\in V$ is an {\sl idempotent} if $a^2=a$.
If $V$ has a unit $\one$, then every
idempotent $a\in V$ gives rise to a vector space direct sum
decomposition of $V$, the {\sl Peirce decomposition}:
$$
V=V_1(a)\oplus V_{1/2}(a)\oplus V_0(a), \qquad 
V_k(a)\colon = \{x\in V \colon ax=kx\}, 
$$
where $k\in \{1,\,\frac{1}{2},\,0\}$ and the corresponding
projectors $E_k(a)\colon V\to V_k(a)$, called {\sl Peirce projectors},
are given by
\begin{equation}\label{id5}
E_1(a)=P(a), \quad E_{1/2}(a)=2L(a)-2P(a), \quad 
E_0(a)=I-2L(a)+P(a).
\end{equation}
If the idempotent $a$ 
is fixed and no confusion is likely
to occur, we write $V_k\colon =V_k(a)$ and 
$P_k\colon =E_k(a)$ for
$k\in \{ 1,\frac{1}{2},0\}$. The Peirce
{\sl multiplication rules} hold
\begin{equation}\label{id6}
\begin{split}
V_0V_0\subset V_0, \qquad V_0V_1&=\{0\}, \qquad
V_1V_1\subset V_1,\\   
(V_0\oplus V_1)V_{1/2}\subset V_{1/2},
&\qquad V_{1/2}V_{1/2}\subset V_0+V_1.
\end{split}
\end{equation}
In particular, $V_0$ and $V_1$ are Jordan subalgebras
of $V$, and  
$[ L(x) \, L(y)]=0$ for $x\in V_0$  and $y\in V_1$ with
the usual commutator product $[\, ,\, ]$. 

A JB-algebra is a {\sl real} Jordan algebra with a
complete norm such that the following conditions hold
$$
\Vert xy\Vert \leq \Vert x\Vert \, \Vert y\Vert , 
\qquad \Vert x^2 \Vert = \Vert x\Vert ^2 , \qquad 
\Vert x^2\Vert \leq \Vert x^2+ y^2\Vert .
$$ 

A JB$^*$-algebra is a {\sl complex} Jordan algebra $U$
with an algebra involution $^* \colon x\mapsto x^*$ and a
complete norm such that the following conditions hold
$$
\Vert xy\Vert \leq \Vert x\Vert \, \Vert y\Vert , 
\qquad \Vert x^*\Vert = \Vert x\Vert , \qquad 
\Vert \{xxx\}\Vert =\Vert x\Vert ^3,
$$
where the {\sl triple product} 
$\{abc\}$ is defined by 
$\{abc\} \colon = (ab^*)c-(ca)b^*+(b^*c)a$ 
and satisfies the {\sl Jordan identity} 
\begin{equation}\label{id4}
\{x\{abc\}y\}= \{\{bax\}cy\}-\{ba\{xcy\}\}+\{\{bay\}cx\}.
\end{equation}
The operators $x\square y\in \mathcal{L}(U)$ are 
defined by $z\mapsto x\square y(z)\colon =\{xyz\}$ for
$z\in U$. 

An element $a\in U$ is {\sl selfadjoint} if
$a^*=a$. The  the set of them, denoted by $U_s$, is a
JB-algebra. Conversely, if $V$ is a JB-algebra then there
is a unique Jordan algebra structure in $U\colon =V\oplus
\imath \, V$, the complexification of $V$, such that
$U_s=V$ and there is a unique norm in $U$ that converts
it into a JB$^*$-algebra \cite{DMW}. We refer to $U$
as the {\sl hermitification} of $V$. The set $\Aut (U)$,
of all Jordan algebra $^*$-automorphisms of $U$, consists
of surjective linear isometries of $U$ and is a
topological group in the topology of uniform convergence
over the unit ball of $U$. By $\Aut ^{\circ}(U)$ we
denote the connected component of the identity in 
$\Aut (U)$. Every element in $\Aut ^{\circ}(U)$ preserves 
the real subspace $V$ and is uniquely determined by this
restriction. 

Let $U$ be a JB$^*$-algebra.  
We write $\Pro (U)$ for the set of self-adjoint idempotents in $U$
and
$\Tri (U)
\colon =
\{a\in U \colon \{a,a,a\}= a\}$ for the set of {\sl
tripotents} in $U$. Clearly
$\Pro (U)\subset \Tri (U)$, and every non zero $a\in \Tri
(U)$ satisfies $\Vert a\Vert =1$. Two elements $a,\,b\in
U$ are {\sl orthogonal} if $ab=0$. A projection $a\in U$ 
is said to be {\sl minimal} if $a\neq 0$ and $P_1(a)U=
\mathbb{C}\,a$, and we let $\Min (U)$ denote the set of
them. For a JB$^*$-algebra it may occur that $\Min
(U)=\emptyset .$

A JBW-algebra is a JB-algebra whose
underlying Banach space $V$ is a dual space, which occurs
if and only if the hermitification $U=V\oplus \imath V$ is
a dual Banach space. In that case $U$ is called a
JBW$^*$-algebra, the predual $U_*$ of 
$U$ is uniquely determined and $\sigma (U,U_*)$, the
weak* topology on 
$U$, is well defined. Let $U$ be the 
JBW$^*$-algebra $U\colon = \mathcal{L}(H)$ of bounded
linear operators $z\colon H\to H$ on a Hilbert space $H$.
Then
$U$ is a unital algebra with plenty of projections each of
which admits a representation of the form 
$$
a= \Sigma_{i\in I}a_i, \qquad \hbox{\rm convergence in
the weak* topology,}
$$
for some indexed family of pairwise orthogonal minimal
projections. The minimal 
cardinal of the set $I$ is the {\sl
rank} of $a$ and $\rank (a)=1$ if and only if $a$ is
minimal. The rank of the algebra $U$ is the rank of its
unit element. 

A JB-algebra $V$ is algebraically (resp. topologically) 
{\sl simple} if $\{0\}$ and $V$ are its only ideals (resp.
closed ideals).

Although not surveyed here, 
we shall occasionally use some relationships between 
JB$^*$-algebras and their associated JB$ ^*$-triples. Our
main reference for JBW-algebras, JBW$^*$-algebras and
JB$^*$-triples are \cite{HAOL} and \cite{UPM}.

\subsection{Manifolds of projections in a JBW-algebra.}

Let $V$ be a JBW-algebra and denote by $U\colon =
V\oplus \imath V$ its hermitification. Then $U$ is a
JBW$^*$-algebra and $\Pro (U)=\Pro (V)$. In the Peirce
decomposition of $U$ induced by 
$a\in \Pro (U)$, the Peirce
spaces are selfadjoint, that is $U_k(a)^*=U_k(a)$, 
and we have $V_k(a)=U_k(a)_s$ where $P_k(a)_{\vert_V}$ is
the Peirce projector of $a$ in the algebra $V$. 
For every $u\in V_{1/2}(a)$, the
operator 
\begin{equation}\label{opg}
G(a,u)\colon =2(u\square a -a\square u)\in
\mathcal{L}(U)
\end{equation}
is an inner derivation of the JBW$^*$-algebra $U$
and the operator-valued mapping
$t\mapsto \exp \, t G(a,u)$, $(t\in \mathbb{R})$,
is a one-parameter group of automorphisms of $U$
each of which preserves $V$. The set $\Pro (V)$,  
endowed with its topology as a subset of $V$, is 
not connected, namely $0$ and $\one $ are 
isolated points. We let $\mathcal{M}(p)$ denote 
the connected component of $p$ in $\Pro (V)$. 
It is
known that $\mathcal{M}(p)$ is a real analytic
manifold whose tangent space at the point $a$ is
$V_{1/2}(a)$, a local chart being given by 
$$
z\mapsto [\exp G(a,z)]\,a, \qquad z\in N,
$$
for $z$ in a suitable neighbourhood $N$ of $0$ in
$V_{1/2}(a)$. As a consequence, all points in 
$\mathcal{M}(p)$ are projections that have the same 
rank as $p$.\begin{mm} As proved in (\cite{KAUJ} th.
4.4),  
$\mathcal{M}(p)$ can also be viewed as a holomorphic 
manifold whose tangent space at the point $a$ is 
$U_{1/2}(a)$ a local chart being given in a 
neighbourhood $M$ of $0\in U_{1/2}(a)$ by 
$$
u\mapsto [\exp 2u\square a]\,a, \qquad u\in M. 
$$ 
For a \end{mm}
projection $p$, the operator $S(p)\colon = \Id -
2P_{1/2}(p)
\in \mathcal{L}(V)$, called the Peirce {\sl
reflection}, is a  {\sl symmetry}, that is, a
selfadjoint involution of
$V$. Namely we have $S(p)x= x$ for $x\in V_1(p)+V_0(p)$ and 
$S(p)x=-x$ for $x\in V_{1/2}(p)$. Besides $S(p)\in \Aut (V)$ 
and $S(p)p=p$, hence $S(p)$ preserves the connected component
$\mathcal{M}(p)$ which in this way is a symmetric manifold. 

We let $\mathfrak{D}$ be the Lie
algebra of all smooth vector fields on
$\mathcal{M}(p)$. A vector field $X$ is now
locally identifiable with a real analytic function
$X\colon N\subset V_{1/2}(a) \to V_{1/2}(a)$. 
We always consider $V_{1/2}(a)$ as submerged in 
$V$. For a function $X\colon \mathcal{M}(p)\to V$, we let
$X_a$ denote the value of $X$ at the point $a\in
\mathcal{M}(p)$. By $X^{\prime}_a$ 
we represent the Fr\'echet derivative of $X$ at
$a$, thus $X^{\prime}_a$ is a continuous linear
operator $V_{1/2}(a)\to V$. For two vector
fields $X,\,Y\in \mathfrak{D}$ we define
\begin{equation}\label{ac}
(\nabla_X\,Y)_a\colon = P_{1/2}(a)\big( Y^{\prime}_a
(X_a)\big), \qquad a\in \mathcal{M}(p).
\end{equation}
It is known that 
$(X, Y)\mapsto \nabla_X\,Y$ is a torsionfree
$\Aut^{\circ}(U)$-invariant affine connection on
the manifold $\mathcal{M}(p)$. For every $a\in 
\mathcal{M}(p)$ and every $u\in V_{1/2}(a)$ the
curve $\gamma_{a,u}(t)\colon = [\exp \,t G(a,u)]\,a$
is a $\nabla$-geodesic that is contained in the
closed real Jordan subalgebra of $V$ generated by
$(a,u)$. Proofs can be
found in \cite{CHIS}.  
 
The following key result is known.
\begin{Theor}\label{frt}
Let $U$ be a 
JBW$^*$-algebra. Then for 
$p\in \Pro (U)$ the 
following conditions are equivalent: 
{\rm (i)} $U_{1/2}(p)$ is a reflexive
space. 
{\rm (ii)} $U_{1/2}(p)$ is linearly
homeomorphic to a Hilbert space.
{\rm (iii)} $\rank U_{1/2}(p)< \infty$. 
For $U\colon =\mathcal{L}(H)$ these conditions are
equivalent to {\rm (iv)}
$\rank (p)<\infty$ or
$\rank (\one -p)<\infty$. 
\end{Theor}
\Dem The equivalence $\rm (i)\iff {\rm (ii)}
\iff {\rm (iii)}$ is known 
\cite{KAUS}. 
The statement concerning $U=\mathcal{L}(H)$ has been
established in \cite{ISTA} as follows:  
From the
expression of the 
Peirce projectors, we have for all $x\in U$ 
\begin{eqnarray*}
P_{1/2}(p)x&=& 2(p\square p -P(p))x=\\
(px+xp)-2pxp&=& p x (\one -p) + (\one -p)xp. 
\end{eqnarray*}
Hence $p\,U_{1/2}(p)\subset U_{1/2}(p)$
and $x\mapsto px$ is a continuous projector 
$U_{1/2}(p) \to p\,U_{1/2}(p)$. Similarly $x\mapsto xp$  
is a continuous projector $U_{1/2}(p)\to U_{1/2}(p)\,p$ 
and since $pU_{1/2}(p)\cap U_{1/2}p=0$ 
we have a topological direct sum decomposition 
$U_{1/2}(p)= X_1\oplus X_2$. Therefore $U_{1/2}(p)$ is 
reflexive if and only if so are the summands. But $X_1\colon
=p\,U_{1/2}(p)= \{ px(\one -p) \colon x\in U\}$ is reflexive
if and only if $\rank (p)<\infty $ or $\rank (\one
-p)<\infty$ as we wanted to see. \qed

Suppose $p\in \Pro (U)$ is such that $\rank
(p)<\infty$. Since this condition is $\Aut
^{\circ}(U)$-invariant, all projections $a\in
\mathcal{M}(p)$ satisfy it too, and the tangent space
$U_{1/2}(a)_s=V_{1/2}(a)$ to
$\mathcal{M}(p)$ at any point $a\in \mathcal{M}(p)$ is a
Hilbert space. In (\cite{DIN}, prop. 9.12) one can
find an explicit expression for an $\Aut
^{\circ}(U)$-invariant scalar product whose norm is equivalent
to $\Vert \cdot \Vert$ in $U_{1/2}(a)$. 
Since $\mathcal{M}(p)$
is connected, that scalar product, denoted
by $\langle \cdot ,\cdot \rangle_a$,  
is determined up to a positive constant coefficient that
can be normalized by requiring that minimal tripotents
should have norm one. We shall not go into details as no
explicit expression of it will be used here. 
\begin{Def}\label{ln}
{\rm We refer to this Hilbert
space norm as the {\sl Levi norm} in $U_{1/2}(a)$ and denote
it  by $\vert \cdot \vert_a$}
\end{Def} 
\begin{Def} {\rm We define a Riemannian metric $g$ on 
$\mathcal{M}(p)$ by} 
$$
g(X, Y)_a\colon = \langle X_a,Y_a\rangle_a, \quad X,\,Y \in
\mathfrak{D}\mathcal{M}(p), \quad a\in \mathcal{M}(p).
$$
\end{Def}\begin{mm}
\begin{Prop}\label{rc}
The affine connection in 
{\rm (\ref{ac})} is the Levi-Civita (respectively,
the K\"ahler) connection on the real
analytic (the holomorphic)
manifold $\mathcal{M}(p)$. 
\end{Prop}
\Dem Indeed, $\nabla$ is compatible with $g$, that is 
$$
X\, g(Y,Z)= g(\nabla_X Y,\, Z)+g(Y,\, \nabla_XZ), \qquad 
X,Y,Z\in \mathfrak{D}\mathcal{M}(p).
$$
Moreover $\nabla$ is torsionfree, hence by (\cite{KLIN}
th.1.8.11) $\nabla$ is the unique Riemann connection on
$\mathcal{M}(p)$. Remark that when $\mathcal{M}(p)$ is 
looked as a holomorphic manifold $\nabla$ is hermitian, that
is, it satisfies $ g(i\,Y,\, i\,Z)= g(Y,\,Z)$, therefore 
$\nabla$ is the only Levi-Civita connection on
$\mathcal{M}(p)$. On the other hand $\nabla_X\,i\,Y=
i\, \nabla_X\,Y$, hence $\nabla$ is the only hermitian
connection on $\mathcal{M}(p)$. Thus the Levi-Civita and the 
hermitian connection are the same in this case and so
$\nabla$ is the K\"ahler connection on $\mathcal{M}(p)$.
\end{mm}


\begin{Theor}\label{str}
The Jordan-Banach algebra $U_0\colon =
\mathcal{L}_0(H)$ of compact operators
and  
$U_2\colon =\mathcal{L}_2(H)$, the
Jordan-Hilbert algebra of Hilbert-Schmidt
operators on $H$, have the same set of
projections and induce on it the same
topology. 
\end{Theor} 
\Dem The first assertion is clear, and $\Pro
(U_0)=\Pro (U_2)$ is precisely the set of
finite rank projections in $\mathcal{L}(H)$. 
It is also clear the above algebras have the
same set of projections of a given rank $r$, 
($1\leq r<\infty$), say $\mathfrak{P}$.
Denote by $\mathfrak{P}_0$ and
$\mathfrak{P}_2$ the Banach manifold
structures defined on $\mathfrak{P}$ according
to our method and to Nomuras' method, 
respectively. The corresponding tangent spaces
at $a$ are 
$$
T_a\mathfrak{P}_0= \{x\in U_0\colon
2ax=x\}, \qquad T_a\mathfrak{P}_2=\{y\in
U_2\colon 2ay=y\}.
$$
Thus we have
$T_a\mathfrak{P}_2\subset T_a\mathfrak{P}_0$.
On the other hand, $a$ itself is a
Hilbert-Schmidt operator and as $U_2$ is an
operator ideal, from $x\in U_0$ and $x=2ax$ we
get $x\in U_2$, hence
$T_a\mathfrak{P}_2= T_a\mathfrak{P}_0$ as
vector spaces. But both
$T_a\mathfrak{P}_2$ and $T_a\mathfrak{P}_0$
are JB$^*$-triples (subtriples of $U_2$ and
$U_0$, respectively) with the same triple
product. It is known that if a Banach space
$X$ admits a JB$^*$-triple structure, then
the triple product determines the topology of
$X$ in a unique way. This in our case implies
that the topologies induced by
$U_2=\mathcal{L}(H)$ and
$U_0=\mathcal{L}_0(H)$ on the tangent space 
$T_a\mathfrak{P}$ coincide. But then also
coincide the topologies induced on
$\mathfrak{P}$ by these two algebras as they
are locally homeomorphic to the same Banach
space. 
\qed


\section{Equations of the geodesics.} 
For any Jordan algebra $V$, we let $V[u,v]$ denote the subalgebra of
$V$ generated by $(u,v)$.  By $\mathbb{S}:=\Sym (\mathbb{R}, 2)$ we
denote the Jordan algebra of the
symmetric 
$2\times 2$ matrices with real entries and the usual Jordan matrix
product. Recall that in $\mathbb{S}$ the set $\Pro (\mathbb{S})$
consists of the isolated points $0$, $\one $, and the one-parameter
family of minimal projections
\begin{equation}\label{ZZ}
B(\theta):= 
\begin{pmatrix}
\cos ^2 \theta & \frac{1}{2}\sin 2 \theta \\
\frac{1}{2}\sin 2 \theta & \sin ^2 \theta
\end{pmatrix}, \qquad \theta \in \mathbb{R}.
\end{equation}
The elements $A:= B(0)$ and $C:= B(\frac{\pi}{2})$ satisfy $A\circ
C=0$ and $A+C= \one$, where $\one $ is the unit of the algebra
$\mathbb{S}$.  The element
\begin{equation}\label{ut}
X:= 
\begin{pmatrix}
0 & 1\\
1&0
\end{pmatrix}\in \mathbb{S}_{1/2}(A)\cap 
\mathbb{S}_{1/2}(C)
\end{equation}
is a non zero tripotent and $\{ A,\, X,\, C\}$ is a basis in $\Sym
(\mathbb{R},\, 2)$.

\begin{Theor}\label{2.1}
  Let $V$ be a unital Jordan algebra and let $a\neq 0$ and $u$ denote,
  respectively, a projection in $V$ and a tripotent in $V_{1/2}(a)$
  such that $au^2=a$.  Then for $V[a, u]$ the following conditions
  hold: {\rm (i)} $V[a, u]= \Lin \,\{a,\, u,\, u^2\}.$ {\rm (ii)} $c:=
  u^2 -a$ is an idempotent such that $ac=0.$ {\rm (iii)} $u^2= a+c$ is
  the unit in $V[a, u].$ {\rm (iv)} There is a unique Jordan
  isomorphism $\psi_{a,u} \colon V[a, u]\to \Sym ( \mathbb{R},\,2)$
  that takes $a,\,u$ and $c$ to $A,\,X$ and $C$ respectively.
\end{Theor}
\Dem From $u\in V_{1/2}(a)$ we get 
$au=\frac{1}{2}u$ and by assumption $au^2=a$.  
Define $c:= u^2 -a$. Then the above
results 
and the fact that $u$ is a
tripotent give 
$$
c^2 = (u^2 -a)^2=u^2 -2au^2+a^2= u^2-a=c.
$$
Therefore $c$ is an idempotent and 
$ac= a(u^2-a)= au^2- a=0$. 
Moreover
$$
cu ^2 =(u^2-a)u^2= u^2- u^2a= u^2-a=c
$$
All this is collected in the following 
table 
$$ 
\begin{tabular}{c|c c c c }
$\circ $ & a & u  & $u^2$ & c \\
\hline
a & a  & $\frac{1}{2}u$ & a & 0\\
u & $\cdot $  & $u^2$ & u & $\frac{1}{2}u$\\
$u ^2$ & $\cdot $ & $\cdot $ & $u ^2$ & c\\
c & $\cdot $  & $\cdot $ & $\cdot $ & c\\
\end{tabular}
$$
which proves that the linear span of the
set $\{a,\, u,\, u^2,\, c\}$ is closed under
the operation of taking Jordan products, and
that $u ^2 = a+c$ is the unit of $V[a,u]$ 
Now it is clear that there is a unique Jordan
isomorphism $\psi_{a,u} \colon V[a,u]\to \Sym
(\mathbb{R},\, 2)$ with the desired
conditions. \qed
\begin{Rem}\label{min}
{\rm For a tripotent $u\in V_{1/2}(a)$, the minimality of 
$a$ is a sufficient (but not necessary) condition for 
$\{auu\}=au^2=a$ to be true.}
\end{Rem}
Indeed, by the Peirce rules $u^2\in V_0(a) \oplus
V_1(a)$ and $V_0(a)\, V_1(a)=\{0\}$, hence 
$au^2\in V_1(a)=\mathbb{R}a$ by the minimality
of $a$, therefore $au^2= \rho a$.
Multiplication by $u$, the fundamental
identities (\ref{id1}) and 
$u^3=u$ yield 
$$
u(au^2) = \rho u a =\frac{\rho}{2}u\qquad 
u(au^2) = (ua) u^2
=\frac{1}{2}u^3=\frac{1}{2}u,
$$
hence $\frac{1}{2}u(\rho -1)=0$ and so $\rho
=1$ since $u\neq 0$. Thus 
$\{auu\}=au^2= a$. \qed 

For pairs $(a,u)$ consisting of a projection
$a\neq 0$ and a tripotent $u\in V_{1/2}(a)$ with 
$au^2=a$ it is quite easy to obtain the equation of 
the geodesic $\gamma_{a,u}(t)$ as shown now. 
Let us define a new product in $V$ via
$x\centerdot y:=
\{xay\}$. Then $(V, \, \centerdot )$ is a unital
Jordan algebra with unit $a$. For $n\in \mathbb{N}$
and $x\in V$ we let $x^{(n)}$ denote the $n$-th power of
$x$ in $(V, \, \centerdot )$. Note that $x^n\neq
x^{(n)}$.  

\begin{Theor}\label{JA} Let $V:= \mathcal{L}(H)$, 
and let $a\neq 0$ 
and $u$ respectively be a projection in $V$ and
a tripotent  in $V_{1/2}(a)$ such that $au^2=a$. If 
$\gamma _{a, u}$  is the geodesic  with
$\gamma (0)=a$ and $\dot \gamma (0)=u$, then
$\gamma _{a, u}(\mathbb{R})\subset V[a,u]$.
More precisely we have 
\begin{equation}\label{eg}
\gamma _{a, u}(t) = (\cos ^2 t )\,a+ (\frac{1}{2} 
\sin 2t)\,u+ (\sin ^2 t )\,u^{(2)},\qquad t\in
\mathbb{R}. 
\end{equation}
\end{Theor}
\Dem Let $G(a,u):= 2(u\square a -
a\square u)\in \mathcal{L}(V)$. We have $\gamma _{a,
u}(t)=[\, \exp \, t G(a,u)\,]\, a$ for all $t\in
\mathbb{R}$,  hence
$\gamma _{a, u}(\mathbb{R})$ is contained in
the closed real linear span of the  sequence $(
G(a,u)^n\,a)_{n\in \mathbb{N}}$. We shall
prove that the assumptions $u^3=u$ and $au^2 =a$  
yield 
$G(a,u)^n\,a\in V[a, u]$ for all
$n\in \mathbb{N}$. 

We have 
$$ G(a,u)\,a=  u\in V[a,u]. $$
By assumption  
$\{auu\}=au^2 =a$, hence  
\begin{eqnarray*}
G(a,u)\,u&=&  
2\{uau\}-2\{ auu\} =2(u^{(2)}-a)\in V[a, u] \\ 
G(a,u)\,u^{(2)}&=& 2\{ u\, a\, u^{(2)}\}- 2 \{ a\, u\,
u^{(2)}\}.
\end{eqnarray*}
By the Peirce multiplication rules
$u^{(3)}= \{u  a u^{(2)}\}\in \{ V_{1/2}(a)\;
V_1(a)\;V_0(a)\} =0$. 
The Jordan identity (\ref{id4}) and (\ref{min}) give 
\begin{eqnarray*}
&{}&\{ a u u^{(2)}\}= \{au\{uau\}\}=
\{u\{auu\}a\}-\{\{uaa\}uu\}+\{ua\{uua\}\}= \\
&{}&\{uaa\}-\frac{1}{2}\{uuu\}+\{uaa\}=\frac{1}{2}u
\in V[a, u]
\end{eqnarray*}
and  
$G(a,u)\, u^{(2)}\in V[a, u].$ 

Note that $a,\, u,\, u^{(2)}$ belong to 
different Peirce $a$-spaces, in particular they are
linearly independent unless $u$ or $u^{(2)}$
vanish. We have assumed $u\neq 0$ and if $\psi$ is
the isomorphism in (\ref{2.1}), then $\psi^{-1}u^{(2)}=
\psi^{-1}\{uau\}=\{XAX\}=C\neq 0.$ Therefore they 
form a basis of $V[a,u]$ and $G(a,u)\,V[a,u]\subset
V[a, u]$. As a consequence
$\gamma  _{a,u}(t)$ has a  unique expression of the
form
$$\gamma _{a, u}(t)= f_1(t) a+
f_{1/2}(t)u+f_0(t)u^{(2)},
\qquad t\in \mathbb{R}, $$ for suitable real analytic
scalar valued functions $f_k(t),\,(k=0,\,1/2,\,1)$. By
taking the  derivative with respect to $t$ and
replacing the expressions previously obtained for 
$G(a,u)\,z$ with $z\in \{ a, \, u, \, u^{(2)}\}$, we
get 
\begin{eqnarray*}
&{}&f^{\prime}_1 (t) a+f^{\prime}_{1/2} (t)
u+f^{\prime}_0 (t) u^{(2)}=  
\dot \gamma _{a, u}(t) = G(a,u)(\gamma_{a, u}(t))=\\
&{}&f_1(t) G(a,u)a +
f_{1/2} (t) G(a,u)u+f_0(t)G(a,u)u^{(2)}=\\
&{}& -2f_{1/2}(t)a +(f_1(t) - 
f_0(t))u+2f_{1/2}(t) u^{(2)},
\end{eqnarray*}
whence we have the
first order ordinary differential equation
$F^{\prime}(t)=A(u)F(t)$  with the initial condition
$F(0)=(1,\,0,\,0)$, where $A(u)$ is a $3\times 3$
constant  (that is, not depending on $t$) matrix and
$F(t)$ is the transpose of $(f_1(t),\,f_{1/2} (t), \,
f_0(t)) $. In fact 
$$
A:=
\begin{pmatrix}
0 &-2 & 0\\
1 & 0 & -1\\
0 & 2 & 0
\end{pmatrix}.
$$
One can easily check that the solution is the curve in
(\ref{eg}). \qed

Motivated by this result, we shall now try to weaken 
the restrictions on $u$. 
\begin{Prop}\label{ppp}
Let $V$, $a$ and $u$ respectively be a unital
JB-algebra, a projection in $V$ and a
tripotent in $V_{1/2}(a)$. Then the following
conditions hold: 
{\rm (i)} $p:= au^2$ is a projection that
satisfies $p\leq a$, $u\in V_{1/2}(p)$ and
$pu^2=p.$ 
{\rm (ii)} For $u\neq 0$ we have $p\neq 0$. 
{\rm (iii)} If $u$ and $v$ are orthogonal as
tripotents in 
$V_{1/2}(a)$, then $p:=au^2$
and $q:=a v^2$ are orthogonal as
projections in $V$ and $pv=uq=0$. 
\end{Prop}
\Dem To see that $p$ is an
idempotent it suffices to consider the
subalgebra of $V$ generated by $a, u$ and the
unit $e$, which is a special algebra.
By the Peirce rules, $u^2\in V_1\oplus V_0$
and so we may
  write $u^2=x+y$ for $x\in V_1$ and $y\in V_0$.  Thus $p=ax+ay = x\in
  V_1$ and $pa=p$.  On the other hand, since $u$ is a tripotent
  $u^2=u^4$ and so using the Peirce rules again
  $$p= au^4=ax^2 + 2a(xy)+ ay^2 = x^2=p^2.$$
  That is, $p$ is a
  projection.  Notice that $(a-p)^2=a-p$ and so $a-p\ge 0$.  From
  $u\in V_{1/2}(a)$ and the fundamental formulas (\ref{id1}) we get
  $up= u(au^2)= (au)u ^2= \frac{1}{2}uu^2= \frac{1}{2}u$ hence $u\in
  V_{1/2}(p)$ and $p\neq 0$ if $u\neq 0$.  To complete (i), we notice
  that \[ pu^2 = p(x+y) = px =p^2 =p\] since $p=x\in V_1(a)$ and $y\in
  V_0(a)$.
  
  Assume now that $u\square v=0$ and set $p:=au^2$, $q:=av^2$.  (We
  recall that two tripotents $u$ and $v$ are orthogonal if and only if
  $\{u,v,v\}=0$ or equivalently $\{u,u,v\}=0$.)  Then $0=\{uve\}=uv$
  and also $u^2v=\{u,u,v\}=0$.  Although $u^2$ is a tripotent and
  $vu^2=0$, we remark that two tripotents $e$ and $f$ which are
  orthogonal in the algebra sense ($ef=0$) may not be orthogonal
  tripotents.  However in our case 
  $$\{u^2,u^2,v\}=u^4v -(vu^2)u^2 + (vu^2)u^2= uv=0$$
  and so $u^2$ and $v$ are orthogonal tripotents.  In particular,
  $u^2v^2=\{u^2,v^2,1\}=0$.

Let
$$
u^2=x+y\in V_1(a)+V_0(a), \qquad 
v^2=x^{\prime} +y^{\prime}\in V_1(a)+V_0(a)
$$
as before. Then $0=u^2v^2= xx^{\prime} +yy^{\prime}$ entails
$xx^{\prime} =yy^{\prime} =0$ and so
$$
pq=(au^2 )(av^2)= [a(x+y)]\,[(x^{\prime} +y^{\prime})a]=
(ax)(x^{\prime} a)=xx^{\prime}=0 .$$
To complete the proof, consider
the product $pv=vp= v(u^2a)= L(v)L(u^2)a$. Since $v$ operator commutes
with $u^2$ and $v\in V_{1/2}(a)$, we have $ pv=
L(v)L(u^2)a=L(u^2)L(v)a= u^2(va)= \frac{1}{2}u^2v=0$ as seen before.
Similarly $uq=0$.  \qed

The following can be considered as a generalization of (\ref{2.1}).
\begin{Theor}\label{sau}
  Let $V:=\mathcal{L}(H)$, and let $a$ and $u$ respectively be a
  finite rank projection in $V$ and a vector in $V_{1/2}(a)$. Then we
  have a finite direct sum decomposition
\begin{equation}\label{ji}
V[a,u]\approx V_0\oplus \big(\bigoplus_1^s V[a_k,u_k]\big)
\end{equation}
where $a_0$ and $a_k$, $(1\leq k\leq s)$, are projections with $a=a_0
+\Sigma a_k$, $u_k$ are tripotents in $V_{1/2}(a_k)$, the subalgebras
$V[a_k,u_k]$ are pairwise orthogonal and $a_0u=0$.
\end{Theor}
\Dem The hermitification $U:= V\oplus iV$ of $V$ is JBW$^*$-triple
and, by the Peirce rules, $U_{1/2}(a)$ is a JBW$^*$-subtriple which
has finite rank by (\ref{frt}). Hence by \cite{KAUS} every element in
$U_{1/2}(a)$ has a unique spectral resolution, that is, a
representation of the form $z= \rho_1 u_1+ \cdots + \rho_su_s$ where
$0<\rho_1< \cdots < \rho_s$, the $u_k$ are pairwise orthogonal
(possibly not minimal) non zero tripotents in $U_{1/2}(a)$ and
$s\leq r= \rank U_{1/2}(a)<\infty$. If $z$ is selfadjoint
(that is,
$z\in V_{1/2}(a)$), then the $u_k$ in the spectral resolution
of $z$ are also selfadjoint. Indeed, $u_j\square u_k=0$ for
$j\neq k$, hence the successive odd powers $z^{2l+1}$ of $z$
are
$$
z^{2l+1}= \rho_1^{2l+1}u_1+\cdots +\rho_su_s^{2l+1}, \qquad (0\leq
l\leq s-1).
$$
and the Vandermonde determinant $\det (\rho_k^{2l+1})$ does not
vanish since the $\rho_k$ are pairwise distinct.  Therefore the $u_k$
are linear combinations with real coefficients of the powers
$z_k^{2l+1}\in V_{1/2}(a)$ 
and so $u_k\in V_{1/2}(a)$.

Now we discuss the algebra $V[a,u]$. Let $u= \xi_1u_1+\cdots
+\xi_su_s$ with $ 0<\xi_1<\cdots <\xi_s$ be a spectral resolution of
$u$, and let $a_k:= au_k^2$ for $1\leq k\leq s$. By (\ref{ppp}) the
projections $a_k$ are pairwise orthogonal
and 
 satisfy
$$
\Sigma a_k\leq a, \quad u_k\in V_{1/2}(a_k),  \quad a_ku_k^2= a_k
\quad (1\leq k\leq s).
$$
Set 
\begin{equation}\label{dpd}
a_0:=a- \Sigma a_k, \qquad V_0:=
\mathbb{R}\,a_0.
\end{equation}
Hence $\dim V_0=1$ at most. We shall see that $ V[a,u]\approx
V_0\oplus \bigoplus V[a_k,u_k]$ as an orthogonal direct sum. For that
purpose consider the successive powers $u^l$ of $u$ which are given by
$ u^l= \xi_1^lu_1^l+\cdots +\xi_s^lu_s^l$, $(1\leq l\leq s)$, since
the $u_k$ are pairwise orthogonal. A Vandermonde argument shows that
$u_k\in V[a,u]$ for $1\leq k\leq s$.  Then $a_k=au_k^2\in V[a,u]$ and
$a_0\in V[a,u]$.  Therefore $ V_0\oplus \bigoplus V[a_k,u_k]\subset
V[a,u]$ On the other hand, from $a=a_0+\Sigma a_k$ and $u=\Sigma
\xi_ku_k$ it follows $V[a,u]\subset V_0\oplus \bigoplus V[a_k,u_k]$
whence we get (\ref{ji}) as soon as we show that the summands satisfy
the required orthogonality properties .

We have already shown that $a_ka_j=0=u_ku_j$ for all $k\neq j$. By
(\ref{ppp}) we have $a_ku_j=0=a_ju_k$ and so the subalgebras
$V[a_k,u_k]$ 
and $V[a_j,u_j]$ are orthogonal for $k\neq j$. It remains
to prove that $a_0u=0$.  By assumption $u\in V_{1/2}(a)$ and from
$u_k\in V_{1/2}(a_k)$ we get $u\in V_{1/2}(\Sigma a_k)$, hence
$$
\frac{1}{2}u=au=(a_0+\Sigma a_k)\,u= a_0u+ (\Sigma a_k)u=
a_0u+\frac{1}{2}u
$$
which completes the proof. \qed
\begin{Cor}\label{egg}
Let $V=\mathcal{L}(H)$, and let $a$ and $u$ be respectively 
a finite rank projection in $V$ and  a vector in
$V_{1/2}(a)$. If $a=a_0+\Sigma a_k$ and 
$u=\Sigma u_k$ are the decompositions given in (\ref{sau}) 
then 
$$
[\exp \,t G(a,u)]\,a= a_0+\Sigma_k [\exp \,t
G(a_k,u_k)]\,a_k
$$
\end{Cor}
\Dem The linearity of $G$ and orthogonality properties of
the elements involved give
\begin{eqnarray*}
G(a,u)&=&G(a_0,u)+\Sigma_k G(a_k,u_k)=\Sigma_k
G(a_k,u_k),\\
&{}& G(a_k,u_k)\, V[a_j,u_j]=0,\qquad (k\neq j).
\end{eqnarray*}
Therefore $G(a,u)^n= \Sigma_kG(a_k,u_k)^n$ for all 
$n\in \mathbb{N}$, and the claim follows from the
definition of exponential mapping. \qed


\section{Geodesics connecting two given points.}

\begin{Prop}\label{en2}
  Let $V$ be a JB-algebra and let $a,b$ be two projections in $V$ with
  $\{aba\}=\lambda a$ and $\{bab\}=\mu b$ for some real numbers
  $\lambda ,\mu$. Then $0\leq \lambda = \mu \leq 1$. Furthermore
  $\lambda =0$ if and only if $ab=0$, and $\lambda =1$ if and only if
  $a=b$.
\end{Prop}
\Dem Projections are positive elements in $V$, hence $\{aba\}\geq 0$
by (\cite{HAOL} prop.  3.3.6). Then $\{aba\}=\lambda a$ entails
$\lambda \geq 0$. But $\lambda \leq 1$ since
$$
\lambda = \Vert \lambda a\Vert =\Vert \{aba\}\Vert \leq \Vert
a\Vert ^2 \, \Vert b\Vert \leq 1.
$$
By (\cite{HAOL} lemma 3.5.2) we have $\Vert \{ab^2a\}\Vert =\Vert
\{ba^2b\}\Vert $ for arbitrary elements $a,b$ in $V$, hence in our
case
$$
\lambda =\Vert \{aba\}\Vert =\Vert \{bab\}\Vert =\mu .
$$
Clearly $\lambda =0$ is equivalent to $aba=0$ which by (\cite{HAOL}
lemma 4.2.2) is equivalent to $ab=0$. For arbitrary projections $p,q$,
the condition $pqp=p$ is equivalent to $p\leq q$, therefore $\lambda
=1$ yields $aba=a$ and $bab=b$ that is $a\leq b$ and $b\leq a$ and so
$a=b$ and conversely. \qed
\begin{Prop}\label{en1} Let $V$ be a unital Jordan 
  algebra and let $a,b\in V$ be two projections such that
  $P(a)b=\lambda a$ and $P(b)a=\lambda b$ hold for some real number
  $0<\lambda < 1$. Then for $V[a,b]$ the following conditions hold:
  {\rm (i)} $V[a, b]= \Lin \,\{a,\, b,\, ab\}$.  {\rm (ii)} There is a
  unique Jordan isomorphism $\phi_{a,b}\colon V[a,b]\to\Sym
  (\mathbb{R},\,2)$ that takes $a,\,b$ and $ab$ respectively to $A,\,
  B(\theta)$ and $A\circ B(\theta)$.
\end{Prop}
\Dem Set $p\colon =ab$. It follows from Macdonalds' theorem that
$$
p^2= \frac{1}{2}a\,\{bab\}+\frac{1}{4}\{ab^2a\}+
\frac{1}{4}\{ba^2b\}
$$
hence in our case $p^2=\frac{\lambda}{4}(a+b+2p)$.  The above
results are shown in the following table
$$
\begin{tabular}{c|c c c}
$\circ $ & a & b  & p\\
\hline
a & a  & p & $\frac{1}{2}(p+\lambda a)$ \\
b & $\cdot $ & b &  $\frac{1}{2}(p+\lambda b)$\\
p & $\cdot $ & $\cdot $ &
$\frac{\lambda}{4}(a+b+2 p)$ \\
\end{tabular}
$$
This shows that the real linear span of the set $\{a,\,b,\,p\}$ is
closed under the operation of taking Jordan products, and so $V[a,b]=
\mathbb{R}a\oplus \mathbb{R}b\oplus \mathbb{R}p.$ It is not difficult
to check that $\{a,\,b,\,p\}$ is a basis for $V[a,b]$ hence $\dim
V[a,b]=3$. The other assertion is now clear. \qed

\begin{Rem}\label{rd}
  {\rm The angle $\theta$ appearing in (\ref{en1}) can be
    expressed in terms of
    $a,\,b$. Indeed, since $\phi_{a,b}$ preserves triple products and
    $\{A,B(\theta),A\}= (\cos^2\theta )A$ we have $\lambda
    =\cos^2\theta $ or $\cos^2\theta
    =\Vert P_1(a)b\Vert$.} 
\end{Rem} 
\begin{Rem}\label{mp} {\rm The conditions $P(a)b=\lambda a$ and
    $P(b)a= \lambda b$ with $0<\lambda <1$ are automatically satisfied
    by any pair of minimal projections $a,b$ with $a\neq b$ and
    $ab\neq 0$.  However, minimality is not necessary in order to have
    them.}
\end{Rem}

We use the isomorphisms $\psi_{a,u} \colon V[a,u]\to\Sym (\mathbb{R},
2)$ and $\phi_{a,b} \colon V[a,b]\to \Sym (\mathbb{R}, 2)$ to show
that two distinct minimal projections $a, \,b$ in $V= \mathcal{L}(H)$
can be joined by a geodesic in $\mathfrak{M}$. \begin{mm}
\begin{Theor}\label{GC} 
  Let $V= \mathcal{L}(H)$ and let $\mathfrak{M}(1)$ be the set of minimal
projections in $V$.  If $a, b$ in
  $\mathfrak{M}(1)$ are such that $a\neq b$, $ab\neq 0$, and $W:=V[a,b]$
  then there exists a tripotent $u\in W_{1/2}(a)$ (unique up to sign)
  such that the geodesic $t\mapsto \gamma_{a,u}(t)$ connects $a$ with
  $b$ in $\mathfrak{M}(1)$. Moreover, $u$ is uniquely
determined by the additional property $b= \gamma_{a,u}(t)$
for some $t>0$. 
\end{Theor}
\Dem The pair of projections $a, \,b$ determines uniquely the algebra
$V[a,b]$ and the Jordan isomorphism $\phi_{a,b}\colon V[a,b]\to \Sym
(\mathbb{R},\,2)$ with the conditions in (\ref{en1}). In particular
$$
\phi (b)=
\begin{pmatrix}
  \cos ^2\theta & \frac{1}{2}\sin 2\theta \\
  \frac{1}{2}\sin 2\theta & \sin ^2 \theta
\end{pmatrix}
$$
for some $\theta$ with $0< \theta <\frac{\pi}{2}$.  Let
$u=\phi_{a,b}^{-1}(X)$ where $X$, given in~\eqref{ut}, is the unique (up
to sign) tripotent in $\mathbb{S}_{1/2}(A)$.  Note that
$au^2=a$.  By (\ref{JA}), $\gamma_{a,u}(\mathbb{R})\subset
V[a,u]\subset V[a,b]$ is a geodesic whose image by the isomorphism $\psi_{a,u} \colon
  V[a,b]\to \Sym (\mathbb{R},2)$ is 
$$\psi(\gamma_{a,u}(t))=
\begin{pmatrix}
  \cos ^2t & \frac{1}{2}\sin 2t \\
  \frac{1}{2}\sin 2t & \sin ^2 t
\end{pmatrix}$$ 
Since $G(a,u)$ is real linear on $u$, the definition of the
exponential function gives
$\gamma_{a,\rho u}(t)=\gamma_{a,u}(t\rho)$ for all $\rho $ and 
$t\in \mathbb{R}$. In particular $\gamma_{a,-u}(t)=\gamma_{a,u}(-t)$. 
A glance to the above matrices shows that either $b=
\gamma_{a,u}(\theta)$ or $b=\gamma_{a,\,-u}(\theta)$ where $\theta >0$. 
\qed 
\begin{Cor}\label{cst} 
    With the notation and conditions in the statement of {\rm (\ref{GC})}, 
    there is a unique vector $v$ in $W_{1/2}(a)$ such that $b= \gamma_{a,v}(1)$ 
    where $v=\theta \,u$ for some tripotent $u$ in $W_{1/2}(a)$ with 
    $au^2=a$ and some $\theta$ with $0<\theta <\frac{\pi}{2}$. 
\end{Cor} 
\Dem As proved before, we have $b= \gamma_{a,u}(\theta)$ for a uniquely 
determined tripotent $u\in W_{1/2}(a)$ with $au^2=a$ and the unique 
real number $\theta$ given by $\cos ^2 \theta =\Vert P_1(a)b\Vert$, 
$0<\theta <\frac{\pi}{2}$. Since $\gamma_{a,\theta u}(1)=\gamma_{a,u}
(\theta)=b$, the vector $v\colon = \theta \,u$ clearly satisfies the 
requirements. \qed 
\begin{Cor} With the above notation, the set $\mathfrak{M}(1)$ is connected.
\end{Cor} 
\Dem Fix any $a\in \mathfrak{M}(1)$. Then $\mathcal{N}_a\colon = \{ b\in 
\mathfrak{M}(1) \colon ab\neq 0\}$ is an open set which is pathwise connected
by (\ref{GC}), hence $\overline{\mathcal{N}_a}$ is also connected. But clearly 
$\overline{\mathcal{N}_a}=\mathfrak{M}(1)$. \qed 
\end{mm}

The set of projections in $\mathcal{L}(H)$ that have rank $r$  
is known to be connected for every fixed $r$, $1\leq r < \infty$, 
\cite{LOOS}. 
In order to extend the preceding results, we let 
$\mathfrak{M}(r)$ denote such a set. Suppose that $V$
is a unital Jordan algebra with unit $\one$ and let $p_1,\cdots , p_n$ be
pairwise orthogonal idempotents with sum
$\one$. Define $V_{i,j}\colon = \{p_iVp_j\}$. Then $V_{i,j}=V_{j,i}$
and the vector space direct sum decomposition, called the joint Peirce
decomposition relative to the family $(p_k)$, holds:
$$
V= \bigoplus_{1\leq i\leq j\leq n}V_{i,j}
$$
Besides we have the following multiplication rules
\begin{equation}\label{pr1}
\begin{split}
  V_{i,j}V_{k,l}=0 \;\; \hbox{\rm if}\; \{i,j\}\cap \{k,l\}=\emptyset
  ,&\qquad V_{i,j}V_{j,k}\subset V_{i,k} \;\; \hbox{\rm
    (pairwise distinct $i$, $j$, $k$)}, \\
  V_{i,j}V_{i,j}&\subset V_{i,i}+V_{j,j}\;\; \hbox{\rm (all $i$,
    $j$)}.
\end{split}
\end{equation}
Furthermore we have
\begin{equation}\label{pr2}
\{ V_{i,j}V_{j,k}V_{k,i} \}\subset V_{i,i}\;\;
\hbox{\rm (all $i$, $j$, $k$)}, \quad 
\{ V_{i,j}V_{j,k}V_{i,j}\}=0 \;\;\hbox{\rm 
(pairwise distinct $i$, $j$, $k$)}.
\end{equation}

Our goal now is to study $V[a,b]$, where $a,\,b$ are projections in
$V$ that have the same finite rank. To simplify the notation we set
$W\colon = V[a,b]$.  By (\cite{NOM2}, lemma 2.5), if $V$ is a
topologically simple Jordan JBW-algebra and $a, \,b$ are two
projections in $V$ with the same finite rank, then $\dim V[a,b]<
\infty$.

Let $e$ be the unit of $W$. Since
$f=(a+b-e)^2$ 
is a positive selfadjoint element, it has a spectral resolution in
$W$. Let it be
\begin{equation}\label{eq0}
(a+b-e)^2= \Sigma_1^ \sigma \lambda_j e_j
\end{equation}
where the number of summands is finite, $\lambda_j\geq 0$ and the
$e_j$ are non zero pairwise orthogonal projections in $W$.  We can
also assume that $\lambda_j$ are pairwise {\sl distinct} and that
$e=\Sigma_1^\sigma e_j$ though the $e_j$ may then fail to be minimal
in $W$. Let
\begin{equation}\label{eq4}
W= \bigoplus_{i,j}W_{ij}, \qquad 
a=\Sigma_{i,j}a_{i,j}, \qquad  b=\Sigma_{i,j}b_{i,j}
\end{equation}
respectively be the Peirce decompositions of $W$, $a$ and $b$ relative
to the complete orthogonal system $(e_j)$. Then (\cite{NOM2}, lemma
2.2) we have $W_{i,j}=\{0\}$ for all $i\neq j$.

We set $W^j\colon =W_{j,j}$, $a_j\colon = a_{j,j}$ and $b_j\colon =
b_{j,j}$, $(1\leq j\leq \sigma)$ to shorten the notation.  The
decompositions in (\ref{eq4}) now read
\begin{equation}\label{eq5}
W= \bigoplus_j W^j, \qquad
a=\Sigma_j a_j, \qquad b= \Sigma_ jb_j
\end{equation} 
Since $W=V[a,b]$ is special, there exists a Jordan $^*$-isomorphism
$\omega \colon W\to \mathfrak{W}\subset \mathfrak{A}^J$ for some
associative algebra $\mathfrak{A}$. Fix any such isomorphism. Then the
elements $A:= \omega(a)$, $B:= \omega(b)$ and $F:=
  \omega(f)$ satisfy
\begin{equation}\label{sji}
A\circ F=AF=ABA=FA=F\circ A.
\end{equation} 
These results have been established in \cite{NOM2} in the context of
JH-algebras, but a careful reading reveals that no essential use of
the scalar product in $V$ has been made.
\begin{Prop}\label{en4}
  In the decomposition in {\rm (\ref{eq5})} we have for $j= 1, 2,
  \cdots , \sigma$: {\rm (i)} The $a_j$ are pairwise orthogonal
  projections and so are the $b_j.$ {\rm (ii)} $W^j= V[a_j,\,b_j]$,
  that is, $W^j$ is the subalgebra generated in $W$ by $a_j,b_j$.
  {\rm (iii)} $P(a_j)b_j=\lambda_ja_j$ and $P(b_j)a_j=\lambda_jb_j$,
  where the $\lambda_j$ are as in {\rm (\ref{eq0})}.
\end{Prop}
\Dem By (\ref{pr1}) $W^jW^j\subset W^j$ and $W^jW^k=\{0\}$ for $j\neq
k$, hence from $a^2 =a$ we get
$$
a^2= (\Sigma a_j)^2= \Sigma a_j^2 , \qquad a= \Sigma a_k
$$
Since the $W^k$ are direct summands in $W$ we get $a_k^2=a_k$ and
similarly $b_k^2=b_k.$

  As $a$, $b$ and $e$ are elements of $\bigoplus V_k$ we have
  $V[a,b]\subset \bigoplus V_k$.  Clearly $V_k\subset W^k$ since $a_k,
  b_k \in W^k$ and so
\[ W = V[a,b] \subset {\textstyle\bigoplus} V_k \subset
{\textstyle\bigoplus} W^k = W.\] Therefore $V_k=W^k$ since the sum is
direct. 

To establish the relations in the last assertion, we set
$$
f_j:= (a_j+b_j-e_j)^2, \qquad (1\leq j\leq \sigma).
$$
and note that $f_j\in W^j$. The orthogonality of the $W^j$ and
(\ref{eq0}) yield
\begin{eqnarray*}
f=(a+b-e)^2&=&(\Sigma a_j+b_j-e_j)^2= \Sigma
(a_j+b_j-e_j)^ 2=\Sigma f_j\\
f=(a+b-e)^2&=& \Sigma \lambda_je_j
\end{eqnarray*}
hence $f_j=\lambda_je_j$, $(1\leq j\leq \sigma)$.  Since $W$ is a
special we can transfer the above relations via the Jordan isomorphism
$\omega \colon W\to \mathfrak{W}\subset \mathfrak{A}^J$. The relations
in (\ref{eq5}) via $\omega$ yield
\begin{eqnarray*}
AF=A\circ F&= &\omega (af)=
\omega (\Sigma \lambda _ka_ke_k)=\omega (\Sigma
\lambda_ka_k)= \Sigma \lambda_kA_k,\\
ABA&=& \omega (aba)=
\omega (\Sigma a_kb_ka_k)=\Sigma A_kB_kA_k
\end{eqnarray*}
which via $\omega ^{-1}$ gives $P(a_k)b_k=\lambda_ka_k$ because the
$W^k$ are direct summands. Similarly we can prove $P(b_k)a_k=\mu_ka_k$
with $\mu_k=\lambda _k$. \qed

\medskip Due to $P(a_k)b_k=\lambda_ka_k$ the spectral values
$\lambda_k$ satisfy $0\leq \lambda_k\leq 1$ and we have the three
following possibilities:

Case I: $\lambda_k=0$. This can not occur for more than one index, say
$k=0$. Then $P(a_0)b_0=0$ and $P(b_0)a_0=0$, hence $a_0$ and $b_0$ are
orthogonal. We shall see below that in this case $\rank a_0=\rank
b_0$, (say $n_0$). Thus $W^0$ is isomorphic to the space of the
diagonal matrices
$$
W\approx \mathbb{R}
\begin{pmatrix}
\begin{array}{c|c}
\one &{}\\
\hline 
{}&0
\end{array}
\end{pmatrix}
\oplus \mathbb{R}
\begin{pmatrix}
\begin{array}{c|c}
0&{}\\
\hline
{}&\one
\end{array}
\end{pmatrix}
$$
with the usual Jordan matrix operations. Here $\one$ is the
$n_0\times n_0$ unit matrix.

Case II: $\lambda_k=1$. This can not occur for more than one index
(say $k=1$). Then $P(a_1)b_1=a_1$ and $P(b_1)a_1=b_1$, which means
that $a_1=b_1$, hence $\rank a_1= \rank b_1$ (say $n_1$) and $W^1$ is
isomorphic to the space of the diagonal matrices
$$
W^1\approx \mathbb{R}\left(
\begin{array}{ccc} 
1&{}&{}\\
{}&\ddots&{}\\
{}&{}&1
\end{array}\right)
$$

Case III: $0<\lambda_k<1$. This may occur for several indices $k$ (the
corresponding $\lambda$ being distinct).  Then proposition (\ref{en1})
applies, hence $W^k$ is Jordan isomorphic to $\Sym (\mathbb{R},\,2)$
via the isomorphism in (\ref{en1}). It is now clear that $a_k$ and
$b_k$ are {\sl minimal} in $W^k$. Since different $W^k$ are
orthogonal, $a_k$ and $b_k$ are also minimal in $W$, that is $\rank
a_k=\rank b_k=1$. Since by assumption $a$ and $b$ had the same rank,
we can now conclude that $\rank a_0=\rank b_0$ as announced earlier.
We can now summarize the discussion in the following
\begin{Theor}\label{str}
  Let $V$ be a topologically simple Jordan JBW-algebra and let $a, \,
  b$ be two projections in $V$ that have the same finite rank. If
  $W^0$, $W^1$ and $W^k$ are the algebras described above, then
  $V[a,\,b]$ is Jordan isomorphic to the finite orthogonal direct sum
\begin{equation}\label{dsd}
V[a, \,b]= W^0 \oplus W^1\oplus \bigoplus_{0\neq k\neq 1}
W^k
\end{equation}
\end{Theor}
Given $a,\, b$ in $\mathfrak{M}(r)$ we show that it is possible to
connect $a$ with $b$ by a geodesic.
\begin{Lemma}\label{en5.2}
  Let the algebra $V$, the projections $a,\, b$ and the decompositions
  $a= \Sigma a_k$ and $b= \Sigma b_k$ be as in (\ref{eq5}). Then
   $P_1(a)b=aba=\Sigma a_kb_ka_k= \Sigma
    P_1^k(a_k)b_k.$
\end{Lemma}
\Dem It follows from the facts that $W$ is special and the $W^k$ are
pairwise orthogonal.

\begin{Lemma}\label{en5.3}
  Let the algebra $V$, the projections $a,\, b$ and the decompositions
  $$
  a=a_0+a_1+\Sigma_{0\neq k\neq 1} a_k, \qquad
  b=b_0+b_1+\Sigma_{0\neq k\neq 1} b_k
  $$
  be given by (\ref{dsd}). If $P_1(a)b$ is invertible in the
  algebra $W_1(a)$, then $a_0=b_0=0$, $a_1=b_1$ and $a_kb_k\neq 0$ for
  $0\neq k\neq 1$.
\end{Lemma}
\Dem It follows directly from the properties of the algebras $W^0$,
$W^1$ and $W^k$ that were established in discussion in (\ref{str}) and
the invertibility of $P_1(a)b$. \qed
\begin{Theor}\label{enG}
  Let $V=\mathcal{L}(H)$ and let $a,\,b$ be two projections in $V$
  with the same finite rank $r$. Assume that $P_1(a)b$ is invertible
  in the unital algebra $V_1(a)$. Then there is a geodesic that joins
  $a$ with $b$ in $\mathfrak{M}(r)$.
\end{Theor}
\Dem We may assume $a\neq b$. Consider the algebra $W\colon =
V[a,\,b]$ and the decompositions in (\ref{dsd}). By (\ref{en5.3}) the
invertibility of $P_1(a)b$ in $W_1(a)$ yields
$$
a_0=b_0=0, \qquad a_1=b_1, \qquad \rank (a_k)=\rank (b_k)=1 \,\;
\hbox{\rm for}\; 0\neq k\neq 1
$$
Thus $W^0=\{0\}$ in our case. Let us define $\gamma_1\colon
\mathbb{R}\to W^1$ to be the constant curve $\gamma_1(t)\colon
=a_1=b_1$, and let $r_1\colon = \rank a_1=\rank b_1$.  Clearly
$\gamma_1$ is a geodesic in the manifold $\mathfrak{M}(r_1)$ 
of the projections in $W^1$ that have fixed rank $r_1$.

For $0\neq k\neq 1$ the projections $a_k$ and $b_k$ are non orthogonal
and minimal in $W^k$. Hence by (\ref{GC}) there is a
geodesic, say $\gamma_k$, that joins $a_k$ with $b_k$ in
$\mathfrak{M}^k(1)$, the manifold of minimal projections in $W^k$.
This curve is of the form
$$
\gamma_k(t)= \gamma_{a_k,u_k}(t) = \big[ \exp \,t
G(a_k,u_k)\big]\,a_k, \qquad t\in \mathbb{R},$$
where $G(a_k,u_k) :=2(a_k\square u_k-u_k\square a_k)$ 
for a tangent vector $u_k\in W^k_{1/2}(a_k)$ 
that is determined by the uniqueness properties established in 
(\ref{cst}). In particular $b_k=\gamma_{a_k,u_k}(1)$. 
We claim that 
$$
\gamma(t) := \gamma_1(t) + \Sigma_{k\neq 1} \gamma_k(t), \qquad
t\in \mathbb{R}
$$
is a geodesic that joins $a$ with $b$ in $\mathfrak{M}(r)$.

By construction we have $\gamma _k(\mathbb{R})\subset W^k$. But these
subalgebras are pairwise orthogonal, hence $\gamma (t)$ is a
projection of rank $r= r_1+\Sigma_{k\neq 1} r_k$ for all $t\in
\mathbb{R}$, that is, $\gamma$ is a curve in $\mathfrak{M}(r)$ and
obviously $\gamma (0)=a$, $\gamma (1)=b$. It remains to show that
$\gamma$ is a geodesic, which amounts to saying that $\gamma$ is of
the form
\begin{equation}\label{eqG}
\gamma (t) = \big[ \exp \, t G(a,
u)\big]\, a, \qquad t\in \mathbb{R},
\end{equation}
for some tangent vector $u\in W_{1/2}(a)$, and it is almost clear that
$u := u_1+\Sigma_{k\neq 1}u_k$ will do.  Indeed, the orthogonality of
the $W^k$ and the expression of the Peirce projectors $P_{1/2}^k(a_k)$
and $P_{1/2}(a)$ for special Jordan algebras easily yield the
inclusions $W_{1/2}^k(a_k)\subset W_{1/2}(a)$ and so
$$u=u_1+\Sigma _{k\neq 1}u_k\in W_{1/2}^1(a_1) \oplus \bigoplus_{k\neq
  1}W^k_{1/2}(a_k)\subset W_{1/2}(a)
$$
Still we have to check that the equality in (\ref{eqG}) holds. To
do this, notice that
$G(W^j, W^k)(W) =\{0\}$ for $j\neq k$, which
is an immediate consequence of (\ref{opg}), the orthogonality
properties of the $W^j$ and $W=\oplus W^j$.  As a consequence
$G(a,u)\,a=\Sigma G(a_j,u_j)\, a_j= \Sigma w_j$ where $w_j\colon
=G(a_j,u_j)\, a_j\in W^j$.  Then
\begin{eqnarray*}
G(a,u)^2\,a&=&G(a,u)\,G(a,u)\,a=G(a,u)\Sigma
w_j=\\
\Sigma_j G(a_j,u_j)\, \Sigma_k w_k&=& 
\Sigma_j G(a_j,u_j)\,w_j=\Sigma_j
G(a_j,u_j)^2 a_j
\end{eqnarray*}
and by induction $G(a,u)^n\,a= \Sigma _jG(a_j,u_j)^n \, a_j$ for all
$n\in \mathbb{N}$, hence
$$
[\exp \, t G(a,u)]\, a= \Sigma_j[ \exp \, t G(a_j u_j)] \, a_j
$$ 
which completes the proof. \qed
\begin{mm}
\begin{Rem}\label{nor}
{\rm The geodesic constructed in (\ref{enG}) satisfies certain normalizing
conditions. Indeed, the pair $(a,\,b)$ determines 
uniquely (up to order) pairs $(a_k,\,b_k)$ via the spectral resolution of 
$(a+b-e)^2$ in the unital algebra $W[a,b]$. In turn, these $(a_k,\,b_k)$ 
determine in a unique way tangent vectors $u_k\in W_{1/2}(a_k)$ such 
that $b_k=\gamma_{a_k,u_k}(1)$ for $1\leq k\leq r$. Finally $u=\Sigma_{k\neq 1}
u_k$. These properties single out the curve $\gamma$ in the class of
geodesics that connect $a$ with $b$.}
\end{Rem}
\end{mm}


\section{Geodesics are minimizing curves.} 
Throughout this section $U$ stands for the algebra $\mathcal{L}(H)$
and $V$ denotes its selfadjoint part. Our next task will be to show
that the geodesic $\gamma_ {a, u}$ joining $ a$ with $ b$ in
$\mathfrak{M}(r)$ is a minimizing curve. That will require some
calculus. Let $a$ be fixed in $\mathfrak{M}(r)$, and let $\vert \cdot
\vert$ denote the Levi norm in $V_{1/2}(a)$ (see \ref{ln}).

\begin{Not}\label{na}
  {\rm We set $\mathcal{N}_a\colon = \{ P_1(a)v \in \mathfrak{M}(r)
    \colon \hbox{\rm $v$ is invertible in $V_1(a)$ } \}$.  Clearly
    $\mathcal{N}_a$ is an open neighbourhood of $a$ in
    $\mathfrak{M}(r)$ and $\{ \mathcal{N}_a \colon a\in
    \mathfrak{M}(r)\}$ is an open cover of $\mathfrak{M}(r)$.
    
    By $B_a\colon = \{x\in V_{1/2}(a)\colon \| x\| <\frac{\pi}{2} \}$
    we denote the open ball in $V_{1/2}(a)$ of radius $\frac{\pi}{2}$
    centered at the origin.  Using the {\sl odd functional calculus}
    for the JB$^*$-triple $V_{1/2} (a)$ (see \cite{KAUF}) one can
    define a mapping $\rho \colon B_a\to V_{1/2}(a)$ by
\begin{equation}\label{les}
u\mapsto \rho (u)\colon = \tan u, \qquad u\in B_a.
\end{equation}
Then $\rho$ is a real analytic diffeomorphism of $B_a$ onto
$V_{1/2}(a)$ whose inverse, also defined by the odd functional
calculus, is
$$
u\mapsto \sigma (u) \colon = \arctan u, \qquad u\in V_{1/2}(a).
$$}
\end{Not}
\begin{Def}
  We define $\Phi _a \colon \mathcal{N}_a\to V_{1/2}(a)$ and $\Psi
  _a\colon V_{1/2}(a)\to V$ by
\begin{eqnarray}
\Phi _a(v)\colon &= & 2\,\big( P_1(a)v\big) ^{-1}
P_{1/2}(a)v,\label{lt}\\ 
\Psi _a(u)\colon &=& [\exp \, G(a, \sigma (u))]\, a,\label{ls}
\end{eqnarray}
\end{Def}
\begin{Lemma}\label{lc} With the above notation, 
  $\Phi_a$ and $\Psi_a$ are real analytic $V$-valued functions.
  Furthermore $\Phi_a (\mathcal{N}_a)\subset V_{1/2}(a)$ and $\Psi_a
    (V_{1/2}(a))\subset \mathcal{N}_a$.
\end{Lemma}
\Dem For $v\in\mathcal{N}_a$, $P_1(a)v$ is invertible in $V_1(a)$.
Hence the mapping $v\mapsto (P_1(a)v)^{-1}$ is well defined and real
analytic in $\mathcal{N}_a$. Clearly $v\mapsto P_{1/2}(a)v$ is real
analytic, hence the product of these two functions, that is, $\Phi_a$
is also real analytic and by the Peirce multiplication rules $\Phi_a
(\mathcal{N}_a) \subset V_{1/2}(a)$.

As said before $u\mapsto \tan (u)$ is a real analytic $V$-valued
function, and so is $u\mapsto G(a, \tan (u))a$
 since $G$ is a
continuous real bilinear mapping. The exponential $u\mapsto \exp \,
G(a, \tan (u))$ is an operator-valued real analytic function, hence by
evaluating at $a$ we get $\Psi_a$, a real analytic function.
Let $u\in 
  V_{1/2}(a)$. We have the decompositions
  $$
  a= a_0+\Sigma_k a_k, \qquad u= \Sigma_k \xi_k u_k,
  $$
  given by (\ref{ji}) with the properties in the statement of
  theorem (\ref{sau}). The odd functional calculus and orthogonality
  gives
  $$
  \arctan u= \Sigma_k \theta_k u_k, \;\; \hbox{\rm where}\;\;
  \theta_k := \arctan \xi_k.
  $$
  Hence $G(a, \sigma (u))= \Sigma_k \theta_k G(a_k,u_k)$.  Again
  using the orthogonality properties $G(a_k,u_k)V[a_j,u_j]=0$ for
  $k\neq j$ we see (recalling the proof of \ref{JA}) that
\begin{equation}\label{psi}
\begin{split}
  \Psi_a(u) &= [\exp \, G(a, \sigma(u))]\, a= \Sigma_k [\exp \,
  \theta_k G(a_k,u_k)]\, a_k = \\
  &{}\Sigma_k (\cos^2 \theta_k)\,a_k + \Sigma _k (\frac{1}{2} \sin
  2\theta_k )u_k+ \Sigma_k (\sin ^2 \theta_k)\, u_k^{(2)}.
\end{split}
\end{equation}
An inspection of this formula shows that $P_1(a)\Psi_a(u)= \Sigma_k
(\cos ^2 \theta_k) a_k$.  Since $\arctan u \in B_a$, we have that
$\max\,\theta_k = \|\arctan u\| < \frac\pi2$.  In particular,
$P_1(a)\Psi_a(u)$ is invertible in $V_1(a)$, that is $\Psi_a(u) \in
\mathcal{N}_a$.  Therefore $\Psi_a(V_{1/2}(a))\subset \mathcal{N}_a$
which completes the proof. \qed

\begin{Prop}\label{in} 
  With the above notation, $\Phi_a\colon \mathcal{N}_a\to V_{1/2}(a)$
  is a real analytic diffeomorphism of $\mathcal{N}_a$ onto
  $V_{1/2}(a)$ whose inverse is $\Psi _a$.
\end{Prop}
\Dem First we show that $\Phi_a$ is invertible in a suitable
neighbourhood $W\subset \mathcal{N}_a$ of $a$.  Let us use the
following self-explanatory notation
$$
\Phi_a(v)=2\,\big(P_1(a)v\big)^ {-1}\,P_{1/2}(a)(v)\colon = 2f(v)
\, g(v).
$$
Note that $f(a)=a$ and $g(a)=P_{1/2}(a)(a)=0$. Thus for $h\in
V_{1/2}(a)$ we have
$$
\Phi_a^{\prime}(a)h= 2\,\big( f^{\prime}(a)h\big) \, g(a) + 2f(a)\,
g^{\prime}(a)h=2a P_{1/2}(a)h=h,
$$
that is $\Phi_a^{\prime}(a)=\Id$ which by the
inverse mapping theorem proves the first claim. Let $v\in V_{1/2}(a)$
and $u\in B_a$ be related by $v=\tan (u)$.  A glance at (\ref{psi})
shows
$$
P_1(a) \Psi_a (u)= \Sigma_k(\cos^2 \theta_k)\, a_k\quad\mathrm{and
  }\quad P_{1/2}(a)\Psi_a(u)=\Sigma_k(\frac{1}{2}\sin 2\theta_k)\,u_k.
$$
Therefore since $u_k\in V_{1/2}(a_k)$,
$$
\Phi_a\Psi_a(u)=2\, \big(
P_1(a)\Psi_a(u)\big)^{-1}\,P_{1/2}(a)\Psi_a(u)= 2\, \Sigma_k (\tan
\theta_k )\, a_ku_k = \Sigma_k \xi_k u_k= u.
$$
Hence $\Phi_a\Psi_a= \Id$.  In particular $\Psi _a$ is the
right-inverse of $\Phi_a$, and the inverse of $\Phi_a$ at least in
$W$. By (\ref{lc}), the mappings $\Phi_a\Psi \colon V_{1/2}(a)\to
V_{1/2}(a)$ and $\Psi_a\Phi_a\colon \mathcal{N}_a\to \mathcal{N}_a$
are well defined and analytic in their respective domains. By
(\ref{enG}) any point in $\mathcal{N}_a$ can be joined with $a$ by a
geodesic that is contained in $\mathcal{N}_a$, hence $\mathcal{N}_a$
is an open connected set.  Since $\Psi_a\Phi_a=\Id$ in $W$, we have
$\Psi_a\Phi_a=\Id$ in $\mathcal{N}_a$ by the identity principle. This
completes the proof. \qed

\begin{Prop}\label{lc3} The family of charts $\{(\mathcal{N}_a,\,
  \Phi_a)\colon a\in \mathfrak{M}(r)\}$ is an atlas which defines the
  manifold $\mathfrak{M}(r)$.
\end{Prop}
\Dem It is easy to check that the above family is a real analytic
atlas whose manifold structure is denoted by
$\mathfrak{M}(r)^{\prime}$. To see that $\mathfrak{M}(r)^{\prime}$ is
the same as $\mathfrak{M}(r)$, recall that
$$
\begin{matrix}
  f\colon U\subset V_{1/2}(a)\to \mathfrak{M}(r) &\qquad &\Phi_a\colon
  \mathfrak{M}(r)'\to V_{1/2}(a) \\
  u\mapsto f(u)\colon = [\exp G(a,u)]\,a &\qquad & v\mapsto
  \Phi_a(v)\colon = 2\, \big( P_1(a)v\big)^{-1}\, P_{1/2}(a)v
\end{matrix}
$$
are local charts of $\mathfrak{M}(r)$ and
$\mathfrak{M}(r)^{\prime}$ at the point $a$. The composite map
$F := \Phi_a\circ f $ can be written in the form
  $$
  F(u)=\Phi_a[f(u)]=2\,
  \big(P_1(a)f(u)\big)^{-1}\,P_{1/2}(a)f(u)=G(u)\, H(u)
  $$
  with self-explanatory notation. Then $G(0)=2a, \; H(0)=
  P_{1/2}(a)a=0, \; H^{\prime}(0)=P_{1/2}(a).$ Therefore, for $h\in
  V_{1/2}(a)$ we have
  $$F^{\prime}(0)h= \big( G^{\prime}(0)\, h\big)
  H(0)+G(0)H^{\prime}(0)h= G(0)H^{\prime}(0)h=2a\, P_{1/2}(a)h=h.
  $$
  Thus $F^{\prime}(0)=\Id$. The remaining part of the proof is
  similar. \qed

We are in the position to prove that geodesics in $\mathfrak{M}(r)$
are minimizing curves. For that we consider $\mathfrak{M}(r)$ as
defined by the atlas $\{ (\mathcal{N}_a, \Phi_a)\colon a\in
\mathfrak{M}(r)\}$.
\begin{Theor}\label{PP} Let $\mathfrak{M}(r)$ be the manifold 
  of projections in $V =\mathcal{L}(H)$ that have fixed finite rank
  $r$.  Let $a\in \mathfrak{M}(r)$ and $\mathcal{N}_a$ be as defined
  in {\rm (\ref{na})} Then for every 
  $b$ in $\mathcal{N}_a$, the geodesic joining $a$ with $b$ is a 
  minimizing curve for the Riemann distance in $\mathcal{N}_a$.
\end{Theor}
\Dem We may assume $a\neq b$. The diffeomorphisms
$B_a\buildrel{\rho}\over
\longrightarrow V_{1/2}(a) \buildrel {\Psi _a}\over \longrightarrow
\mathcal{N}_a$ give a unique pair $(u, v)\in B_a\times V_{1/2}(a)$
such that
$$
v=\tan (u) \qquad \Psi _a(v)=b.
$$
There is a unique normalized geodesic $\gamma_{a,u}$ that joins $a$
with
$b$ in the manifold $\mathfrak{M}(r)$ and has initial velocity $u=\dot
\gamma_{a,u}(0)\in B_a$.  In particular we have
$$
b=\Psi_a(\rho (u))= \gamma_{a,u}(1)= [\exp \,G(a,u) ]\,a
$$
and the exponential mapping $\exp \colon B_a\to \mathfrak{M}(r)$ is
a homeomorphism of $B_a$ onto the open set $\mathcal{N}_a$ in
$\mathfrak{M}(r)$. This will allow us to apply the Gauss lemma
(\cite{KLIN} 1.9). For that purpose, we show that $\gamma_{a,u}(t)$
belongs to $\mathcal{N}_a$ for all $t\in [0,1]$. Indeed, the segment
$[0,1]\,u$ is contained in $B_a$, hence its image by $\Psi_a\circ
\rho$ lies in the set $\mathcal{N}_a$. We shall now see that
$$
\Psi_a[\rho (tu)] =\gamma_{a, u}(t), \qquad t\in [0,1].
$$
Let $t\in [0,1]$ and set $v_t\colon =\tan (tu)$. The odd functional
calculus gives
$$
\Psi_a(v_t)= [\exp \, G(a, tu)]\,a= [\exp \,t\, G(a, u)]\,a=
\gamma_{a,u}(t)
$$
as we wanted to see. For the Riemann connection $\nabla$ in
$\mathfrak{M}(r)$, the radial geodesics are minimizing curves (by the
Gauss lemma). Hence it suffices to see that $\gamma_{a,u}$ is a radial
geodesic, which is a consequence of the fact $\gamma_{a,u}[0,1]\subset
\mathcal{N}_a$. \qed

It is now reasonable to ask what can be said about the neighbourhood
$\mathcal{N}_a$. We refer to $\mathcal{O}_a\colon = \{b\in
\mathfrak{M}(r)\colon P_1(a)b \;\,\hbox{\rm is not invertible in} \;\,
V_1(a)\}$ as the {\sl antipodal} set of $a$. Clearly $\mathcal{O}_a$
is a closed subset of $\mathfrak{M}(r)$. 
\begin{Prop} \label{ei} Let $\mathfrak{M}(r)$ be 
  the manifold of all projections in 
$V=\mathcal{L}(H)$ that have a given finite rank $r$. Then for
  any $a\in \mathfrak{M}(r)$ the antipodal set of $a$ has empty
  interior.
\end{Prop} 
\Dem Let $a\in \mathfrak{M}(r)$ and set $K\colon = a(H)\subset H$.
Note that $\hbox{\rm dim}\,K=\hbox{\rm rank}\, a=r<\infty$. The
operators in $V_1(a)= aVa$ can be viewed as operators in $\mathcal{
  L}(K)$, therefore the {\it determinant} function is defined in
$V_1(a)$ and an element $z\in V_1(a)$ is invertible if and only if
$\hbox{\rm det}\, (z)\neq 0$.  The function $b\mapsto \det (P_1(a)b)$
is real analytic on $\mathfrak{M}(r)$. If $\mathcal{O}_a$ has non
empty interior, then $\det (P_1(a)b)$ vanishes in a non void open
subset of $\mathfrak{M}(r)$, which is connected, therefore by the
identity principle the determinant function would be identically null
which is a contradiction.  \qed 

\begin{mm}
We let $\mathfrak{S}(r)$
denote the subgroup of
$\Aut (U)$ generated by the set of Peirce
reflections $\sigma_p$, where $p\in
\mathfrak{M}(r)$. Each $\sigma_p$
preserves $\mathfrak{M}(r)$ and induces a real
analytic symmetry of this manifold. 
\begin{Prop} Let $V=\mathcal{L}(H)$. Then $\mathfrak{M}(r)$ is
homogeneous under the action of any of the groups $\Aut ^{\circ}(V)$
and 
$\mathfrak{S}(r)$.   
\end{Prop}
\Dem Let $a,\,b$ be any pair of points in $\mathfrak{M}(r)$ with
$a\neq b$. If $b\in \mathcal{N}_a$ then by
  (\ref{enG}) there is a geodesic $\gamma_{a,u}(t)= [\exp \, t
G(a,u)]\,a$ that joins $a$ with $b$ in $\mathfrak{M}(r)$. But
$g(t)\colon =[\exp \, t G(a,u)]$ is an element of $\Aut ^{\circ}(V)$
for all $t\in \mathbb{R}$.
   Now consider the case 
   $b\notin \mathcal{N}_a$. By (\ref{ei}) the antipodal set
  of $a$ has empty interior, hence $W\cap \mathcal{N}_a\neq \emptyset$
  for every neighbourhood $W$ of $b$ in $\mathfrak{M}(r)$. If $c\in
  W\cap \mathcal{N}_a$ then we can connect $a$ with $c$ and $c$ with
  $b$ by geodesics in $\mathfrak{M}(r)$ hence we connect $a$ with $b$
  by
  a curve in $\mathfrak{M}(r)$. This completes the proof for  $\Aut
 ^{\circ}(V)$.

Consider a pair of points $a,\, b$
in $\mathfrak{M}(r)$ such that $b\in
\mathcal{N}_a$. Then it is easy to find a
symmetry $\sigma_p$ that exchanges $a$ with
$b$. Namely, there is a geodesic that connects
$a$ with $b$ in $\mathfrak{M}(r)$, say
$\gamma\colon t \mapsto \gamma(t)$, $t\in
\mathbb{R}$. Thus $a=\gamma(0)$ and
$b=\gamma(1)$. Define the {\sl geodesic middle
point} of the pair $(a,\,b)$ as $c\colon 
=\gamma(\frac{1}{2})$, and let $\sigma_c$ be the
Peirce reflection with center at $c$. Clearly
$\sigma_c$ preserves the curve $\gamma$ and
exchanges $a$ with $b$. 

Now let $a,\,b$ be arbitrary in
$\mathfrak{M}(r)$. Since $\mathfrak{M}(r)$ is
a connected locally path-wise connected
topological space, it is globally path-wise
connected. Thus there is a curve $\Gamma\colon t\mapsto \Gamma
(t)$, $t\in [0,1]$, that joins $a$ with $b$ in
$\mathfrak{M}(r)$, and a standard compactness
argument shows that there is a finite sequence
of points $\{x_1,\cdots , x_r\}$ in $\Gamma$
such that $x_1=a$, $x_r=b$ and each
consecutive pair of the $x_k$ can be joined
by a geodesic in $\mathfrak{M}(r)$. For each
pair $(x_k, x_{k+1})$ consider the
geodesic middle point $c_k$ and the
corresponding symmetry exchanging $x_k$ with
$x_{k+1}$. Then the composite
$\sigma_1\circ \cdots \circ \sigma_k$ lies in 
$\mathfrak{S}(r)$ and exchanges $a$ with $b$. 
\qed
\begin{Cor}
Let $a, \,b\in \mathfrak{M}(r)$. Then $b\in \mathcal{O}_a$ 
if and only if $a\in \mathcal{O}_b$. 
\end{Cor}
\Dem Let $\sigma$ be a symmetry that exchanges $a$ with $b$. 
The relation $b\in \mathcal{O}_a$ is equivalent to 
$P_1(a)b$ is not invertible in $V_1(a)$, which applying $\sigma$ is 
converted into $P_1(b)a$ is not invertible in $V_1(b)$ that is 
$a\in \mathcal{O}_b$. 
\qed

\end{mm}
We are now in the position to compute the Riemann distance in $
\mathfrak{M}(r)$.
\begin{Theor}\label{RD} Let $\mathfrak{M}(r)$ be 
  the manifold of all projections in 
  $V=\mathcal{L}(H)$ that have a given finite rank $r$. If $a, \, b$
are points in
  $\mathfrak{M}(r)$ and $\gamma_{a,u}(t)$ is the normalized 
geodesic connecting $a$ with $b$ in $\mathfrak{M}(r)$ then the
Riemann distance between them is
  $$
  d(a,\, b)=\big( \Sigma_1^r
\theta_k^2\big)^{\frac{1}{2}}
 $$
	where $\theta_k=
\cos ^{-1} \big(\Vert P_1(a_k)b_k\Vert
^{\frac{1}{2}}\big)$ and $\Vert\cdot\Vert$ stands for the usual operator 
norm. 
\end{Theor}
\Dem Discard the trivial case $ a= b$. Consider first the case $
b\in \mathcal{N}_a$. By (\ref{enG}, \ref{nor}) we have
$b=\gamma_{a,u}(1)$ for some tangent vector $u=\Sigma u_k$ where $u_k=
\theta_k\,v_k$, $0<\theta_k<\frac{\pi}{2}$, $(1\leq k\leq r)$, and
the $v_k$ are pairwise orthogonal (in the JB$^*$-triple sense)
minimal tripotents in $W_{1/2}(a_k)$. Hence the $v_k$ are
pairwise orthogonal in the Levi sense (see \cite{DIN} prop. 9.12 and
9.13) and so, if $\vert \cdot\vert_a$ denotes the Levi norm in 
the tangent space $V_{1/2}(a)$, we have 
$ \vert u\vert_a ^2 = \Sigma_1^r \theta _k^2$ (recall
that minimal tripotents satisfy $\Vert v\Vert =\vert v\vert_a =1$).
Therefore since the Levi form is
$\Aut ^{\circ}(V)$-invariant and
$\Aut ^{\circ}(V)$ is transitive in $\mathfrak{M}(r)$, 
$$
\vert \dot \gamma _{ a,u}(t)\vert_{\gamma _{ a,u}(t)} = \vert \dot
\gamma _{ a,u}(0)\vert_a= \vert u\vert_a =  
\big( \Sigma_1^r \theta _k^2\big)^{\frac{1}{2}}
$$
for all $t\in \mathbb{R}$, hence
$$ 
d(a, b)= \int_0^1 \vert \dot \gamma_{a,u} (t)\vert
_{\gamma_{a,u} (t)}d t= 
\vert u\vert_a = \big( \Sigma_1^r \theta
_k^2\big)^{\frac{1}{2}} .$$
Here $\theta_k= \cos^{-1}\big(\Vert P_1(a_k)b_k\Vert \big) ^{\frac{1}{2}}$ 
and $\Vert \cdot \Vert$ is the JB$^*$-triple (that is, the usual
operator)  norm. 
In the case $b\in \mathcal{O}_a$, consider a sequence $(b_j)_{j\in
\mathcal{N}_a}$  with $b_j \to b$ in $\mathfrak{M}(r)$ and $a 
b_j\neq 0$ for all $j$.  Applying the above 
to each $j$ and taking the limit we get 
the result. \qed


\end{document}